\documentclass[reqno,12pt]{amsart}

\usepackage{epsf}
\usepackage{graphics}
\usepackage{amssymb}
\usepackage{amsmath}

\date{}

\theoremstyle{plain}
\newtheorem{theorem}{Theorem}

\newtheorem{lemma}{Lemma}

\theoremstyle{definition}

\theoremstyle{remark}

\newtheorem*{examples}{Examples}

\newtheorem*{remark}{Remark}

\def\C{{\mathbb C}}

\def\N{{\mathbb N}}

\title{Unknotting Sequences for Torus Knots} 

\author{Sebastian Baader}

\begin{document}

\begin{abstract} The unknotting number of a knot is bounded from below by its slice genus. It is a well-known fact that the genera and unknotting numbers of torus knots coincide. In this note we characterize quasipositive knots for which the genus bound is sharp: the slice genus of a quasipositive knot equals its unknotting number, if and only if the given knot appears in an unknotting sequence of a torus knot.
\end{abstract}

\maketitle

\section{Introduction}

The unknotting number is a classical measure of complexity for knots. It is defined as the minimal number of crossing changes needed to transform a given knot into the trivial knot \cite{We}. An estimate of the unknotting number is usually difficult, even for special classes of knots such as torus knots. It is the content of the Milnor conjecture, proved by Kronheimer and Mrowka~\cite{KM}, that the unknotting number of torus knots coincides with another classical invariant of knots, the slice genus. The slice genus $g_*(K)$ of a knot $K$ is the minimal genus among all surfaces smoothly embedded in the 4-ball with boundary $K$. In general, the slice genus of a knot $K$ provides a lower bound for its unknotting number $u(K)$:
$$u(K) \geq g_*(K).$$ 
Loosly speaking, this genus estimate tends to be good for quasipositive knots. A knot $K$ is quasipositive, if there exists a quasipositive braid, i.e. a braid which is a finite product of conjugates of the positive standard braid generators $\sigma_i$, whose closure is $K$. Quasipositive knots, or rather links, are precisely the knots or links that arise from intersections of complex plane curves with the unit $3$-sphere in $\C^2$ (\cite{BO}, \cite{Ru1}). In this note we characterize quasipositive knots for which the genus estimate is sharp. For this purpose we introduce the following terminology: an unknotting sequence for a knot $K$ is a finite sequence of knots
$$K=K_n, K_{n-1}, K_{n-2}, \ldots, K_1, K_0,$$
such that:
\begin{enumerate} 
\item $u(K_i)=i$, $0 \leq i \leq n$ (in particular, $K_0$ is the trivial knot), \\
\item two succeeding knots of the sequence are related by one crossing change. \\
\end{enumerate} 
We say that a knot $K$ can be unknotted via a knot $L$, if $L$ is contained in an unknotting sequence for $K$. It should be mentioned that unknotting sequences are highly non-unique: every knot of unknotting number at least two can be unknotted via infinitely many different knots \cite{Ba2}.

\begin{theorem} Let $K$ be a quasipositive knot. Then the equality $u(K)=g_*(K)$ holds, if and only if there exists a torus knot $T(p,q)$ that can be unknotted via $K$. 
\end{theorem}

\begin{examples} The three knots $7_2$, $7_3$ and $7_5$ (in Rolfsen's notation \cite{Ro}) are positive, hence quasipositive by Nakamura and Rudolph  (\cite{Na}, \cite{Ru3}). Moreover, they all satisfy the equality $u=g_*$. According to Theorem~1, they are to be contained in some unknotting sequences for torus knots. Actually, they are all contained in the following two unknotting sequences for the torus knot $T(7,2)=7_1$:

$7_1$, $7_5$, $7_2$, unknot,

$7_1$, $7_3$, $7_2$, unknot.

\noindent
These two sequences are depicted in Figure~1 (the first one on the left; the second one on the right). The knot diagram at the top of the figure is a non-minimal two-bridge diagram for the torus knot $7_1$.

Likewise, every positive twist knot with an odd number of crossings $n$ is contained in an unknotting sequence of the torus knot $T(n,2)$. In contrast, the slice quasipositive knot $8_{20}$ is not contained in any unknotting sequence of a torus knot, since $g_*(8_{20})=0$.
\end{examples}

\begin{figure}[ht]
\scalebox{0.9}{\raisebox{-0pt}{$\vcenter{\hbox{\epsffile{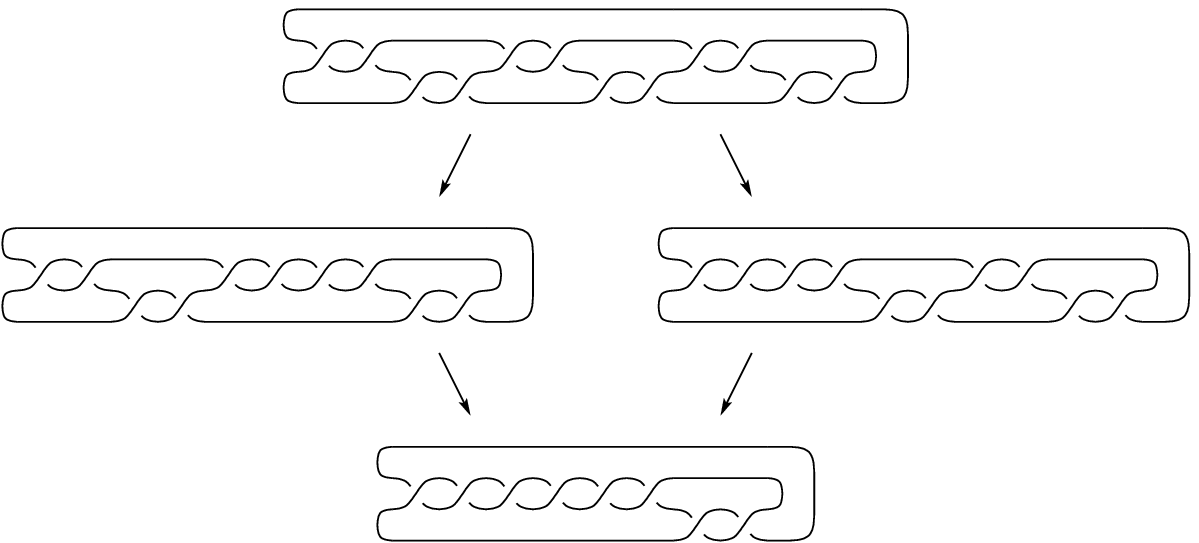}}}$}} 
\caption{}
\end{figure}

\begin{remark} Theorem~1 stays true if we replace the slice genus $g_*$ by any invariant that coincides with $g_*$ on quasipositive knots, for example the concordance invariants $\frac{s}{2}$ or $\tau$ coming from knot Khovanov homology and Floer homology, respectively (see~\cite{HO} for a comparison of these two invariants). 
\end{remark}

The assumption of Theorem~1 that $K$ be quasipositive is essential, as we can see by looking at the figure-8 knot $4_1$ with $u(4_1)=g_*(4_1)=1$: if the figure-8 knot were contained in the unknotting sequence of a torus knot, then its unknotting number would have to coincide with its Rasmussen invariant $s(4_1)=0$. Indeed, the concordance invariant~$s$ gives a sharp bound for the unknotting number of torus knots, hence for all knots of their unknotting sequences~\cite{Ra}.

There do, however, exist many non-quasipositive knots that are contained in unknotting sequences of torus knots. For example, the torus knot $T(5,2)=5_1$ can be unknotted via the knots $8_2$, $8_7$. According to the classification of small quasipositive knots \cite{Ba1}, these are both non-quasipositive knots.

The `if'-part of Theorem~1 is obviously true: the equality $u=g_*$ holds for all torus knots. Further, a single crossing change cannot decrease the slice genus $g_*$ by more than one. Therefore, the equality $u=g_*$ holds for all knots of their unknotting sequences. The `only if'-part is a direct consequence of the following two lemmas, which we shall prove in Sections~2 and~3, respectively.

\begin{lemma} Let $K$ be a quasipositive knot. If the equality $u(K)=g_*(K)$ holds, then there exists a positive braid knot that can be unknotted via~$K$. 
\end{lemma}

Here a positive braid knot is the closure of a finite product of positive standard braid generators $\sigma_i$.

\begin{lemma} For every positive braid knot $K$ there exists a torus knot $T(p,q)$ that can be unknotted via $K$.
\end{lemma}

Concerning Lemma~2, we should mention that the equality $u=g_*$ holds for all positive braid knots. This follows easily from results of Kawamura \cite{Ka} and Boileau-Weber \cite{BW}.

\section{From quasipositive knots to positive braids}

One fundamental consequence of Kronheimer and Mrowka's result is the so-called slice-Bennequin inequality~\cite{Ru2}:
\begin{equation}
g_*(K) \geq \frac{1}{2}(1+w(D)-n(D)).
\label{ineq1}
\end{equation}
Here $K$ denotes the closure of a braid diagram $D$ with $n(D)$ strands and algebraic crossing number $w(D)$. The latter is the number of positive generators minus the number of negative generators of the braid $D$. For quasipositive braid diagrams, this inequality is actually an equality~\cite{Ru2}.

\begin{proof}[Proof of Lemma~1]
Let $K$ be a quasipositive knot with $u(K)=g_*(K)$. By the above observation, $K$ has a quasipositive braid diagram $D$ with 
$$u(K)=g_*(K)=\frac{1}{2}(1+w(D)-n(D)).$$
Changing one negative crossing of $D$ into a positive one, we obtain a diagram $D'$ with $w(D')=w(D)+2$. The unknotting number of the corresponding knot $K'$ increases by one, by (\ref{ineq1}) and the genus estimate. Therefore $K'$ can be unknotted via $K$. Continuing in this manner, we end up with a positive braid knot which can be unknotted via $K$.
\end{proof}

\section{From positive braids to torus knots}

As we already mentioned, the unknotting number and the slice genus are equal for any knot $K$ with a positive braid diagram $D$. Further, these numbers coincide with Bennequin's bound:
\begin{equation}
u(K)=g_*(K)=\frac{1}{2}(1+w(D)-n(D)).
\label{ineq2}
\end{equation}
In this case $w(D)$ is simply the number of crossings of the diagram $D$, since these are all positive. The following observation is an easy consequence of (\ref{ineq2}): if we apply one of the two local moves shown in Figure~2 to a positive braid knot $K$, then we obtain another positive braid knot $K'$ that can be unknotted via $K$. Indeed, the unknotting number increases precisely by half the number of positive crossings (or equivalently, by the number of crossing changes) we insert by the moves, by (\ref{ineq2}). We call the upper and lower move a right turn and left turn, respectively.

\begin{figure}[ht]
\scalebox{0.9}{\raisebox{-0pt}{$\vcenter{\hbox{\epsffile{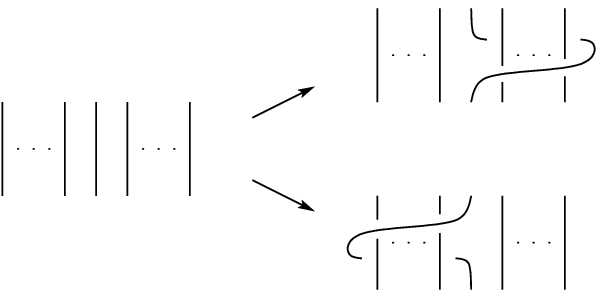}}}$}} 
\caption{}
\end{figure}

\begin{proof}[Proof of Lemma~2]
Let $K$ be a knot with an $n$ strand positive braid diagram $D$. We will see that for $k \in \N$ large enough, the knot $K$ is contained in an unknotting sequence of the torus knot $T(n,kn+1)$. For this purpose, we introduce special braid diagrams for torus knots of type $T(n,kn+1)$, where the $k$ full twists of the individual strands are separated, as shown in Figure~3 for the case $T(4,13)=T(4,3 \cdot 4+1)$. There the three full twists account for the `$3 \cdot 4$', whereas the three crossings at the bottom account for the `$+1$'.

\begin{figure}[ht]
\scalebox{0.9}{\raisebox{-0pt}{$\vcenter{\hbox{\epsffile{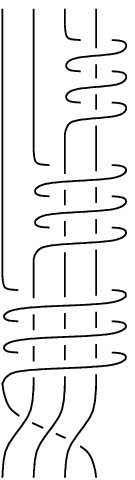}}}$}} 
\caption{}
\end{figure}

The procedure that we will describe now is illustrated in Figure~5. The three arrows therein correspond to the individual steps of the procedure.

Let $s_1$ be the first strand of the braid diagram $D$, i.e. the strand starting at the bottom left of $D$ and ending at a certain position $\tau_1 \geq 2$. Applying a suitable number of left and right turns to $s_1$, we obtain another positive braid of the form:
$$(\sigma_1 \sigma_2 \cdots \sigma_{n-1} \sigma_{n-1} \sigma_{n-2} \cdots \sigma_1)^{k_1} \beta_1 \sigma_1 \sigma_2 \cdots \sigma_{\tau_1-1},$$
for some $k_1 \in \N$ and some positive braid word $\beta_1$ not involving the first generator $\sigma_1$. In fact, $\beta_1$ is the braid that arises from the original braid by deleting the first strand. A conjugation with the braid word $\sigma_1 \sigma_2 \cdots \sigma_{\tau_1-1}$ transforms the above braid into another one representing the  same knot:
$$ \sigma_1 \sigma_2 \cdots \sigma_{\tau_1-1} (\sigma_1 \sigma_2 \cdots \sigma_{n-1} \sigma_{n-1} \sigma_{n-2} \cdots \sigma_1)^{k_1} \beta_1.$$
The latter in turn is isotopic to the following positive braid:
$$ \sigma_1 (\sigma_1 \sigma_2 \cdots \sigma_{n-1} \sigma_{n-1} \sigma_{n-2} \cdots \sigma_1)^{k_1} \sigma_2 \cdots \sigma_{\tau_1-1} \beta_1$$
(this isotopy is shown in Figure~4; the whole step is illustrated by the first arrow of Figure~5).

\begin{figure}[ht]
\scalebox{0.9}{\raisebox{-0pt}{$\vcenter{\hbox{\epsffile{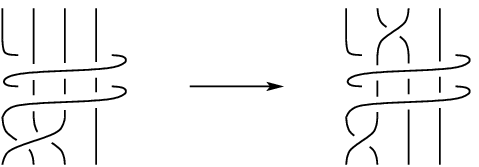}}}$}} 
\caption{}
\end{figure}

Now we restart the whole procedure, let the first strand of the upper part of the braid $\sigma_2 \cdots \sigma_{\tau_1-1} \beta_1$ wind a few times around all the other strands except $s_1$ and obtain a positive braid of the form:
$$ \sigma_1 (\sigma_1 \cdots \sigma_{n-1} \sigma_{n-1} \cdots \sigma_1)^{k_1}
(\sigma_2 \cdots \sigma_{n-1} \sigma_{n-1} \cdots \sigma_2)^{k_2} \beta_2 
\sigma_2 \sigma_3 \cdots \sigma_{\tau_2-1},$$
for some $k_2, \tau_2 \in \N$, $\tau_2 \geq 3$, and some positive braid word $\beta_2$ not involving the first two generators $\sigma_1$, $\sigma_2$.
Again, we may replace the above braid by another one representing the same knot:
$$ \sigma_2 \sigma_1 (\sigma_1 \cdots \sigma_{n-1} \sigma_{n-1} \cdots \sigma_1)^{k_1}
(\sigma_2 \cdots \sigma_{n-1} \sigma_{n-1} \cdots \sigma_2)^{k_2}
\sigma_3 \cdots \sigma_{\tau_2-1} \beta_2$$
(see the second arrow of Figure~5). 

Iterating this procedure until we reach the last strand of the braid, we obtain a positive braid diagram which is almost a diagram of a torus knot of type $T(n,kn+1)$, except that the number of full twists $k_i$ of the individual strands may vary from strand to strand. However, they can easily be made equal by adding a few more twists, where necessary.
\end{proof}

\begin{figure}[ht]
\scalebox{0.8}{\raisebox{-0pt}{$\vcenter{\hbox{\epsffile{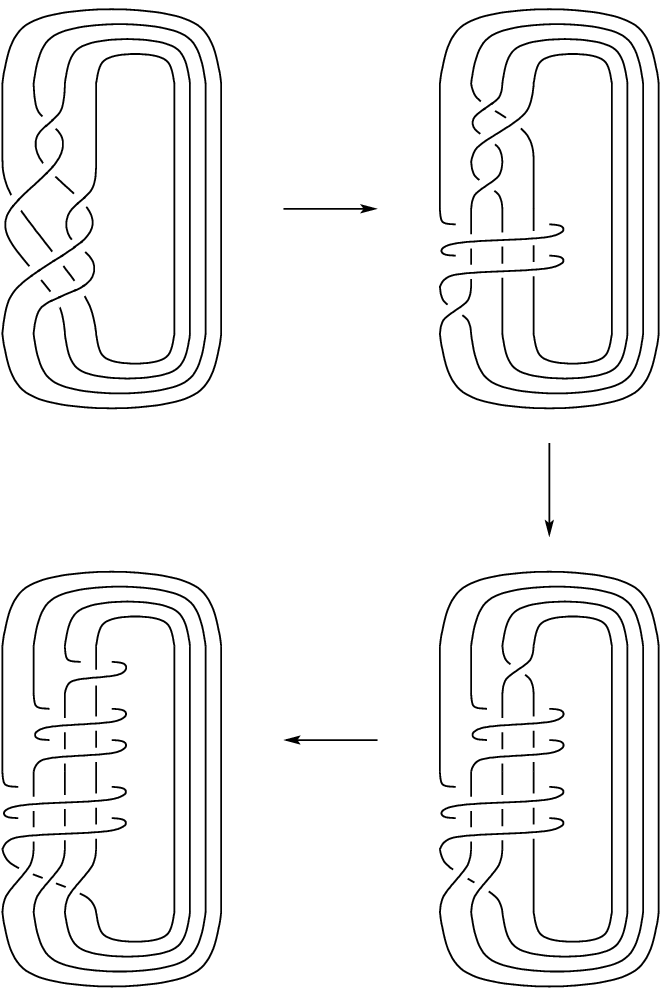}}}$}} 
\caption{}
\end{figure}

\bigskip
\noindent
Department of Mathematics,
ETH Z\"urich, 
Switzerland

\bigskip
\noindent
\emph{sebastian.baader@math.ethz.ch}

\end{document}